\documentclass[11pt]{article}
\usepackage{amsmath,amssymb,amsfonts}
\usepackage{graphicx}
\newcommand{\Z}{\mathbb{Z}}
\newcommand{\R}{\mathbb{R}}

\newcommand{\lo}{\longrightarrow}

\newcommand{\fg}{\mathfrak g}

\newtheorem{definition}{Definition}[section]
\newtheorem{lemma}[definition]{Lemma}
\newtheorem{theorem}[definition]{Theorem}

\newtheorem{cor}[definition]{Corollary}

\newtheorem{pro}[definition]{Proposition}

\def\Af {Aff_\circ(\R)}
\def\co {cohomogeneity one }

\begin{document}
\begin{center} {\bf \Large Cohomogeneity one three dimensional anti de Sitter space,\\
 proper and nonproper actions}\\
{\bf P. Ahmadi  }\\
{ \it Department of Mathematics, Faculty of Sciences, University of Zanjan \\
P.O.Box 45195-313 Zanjan-Iran }\\
email: p.ahmadi@znu.ac.ir

\end{center}
{\small \leftskip 2cm \rightskip 2cm \baselineskip 5mm
\begin{abstract}
In this paper we give a classification of closed and connected Lie
groups, up to conjugacy in $Iso({\bf adS_3})$, acting by \co on the three
dimensional anti de sitter space ${\bf adS_3}$. Then we determine the
causal character of the orbits and the orbit spaces, up to homeomorphism, in both cases, proper and nonproper
actions. When the action is proper, we show that there is no exceptional orbit and causal characters of the principal orbits are the same.
\end{abstract}
\vspace{5mm} {\bf Keywords}:  {\it Cohomogeneity one , Anti de
Sitter space.}}
 \vspace*{5mm}
\section{Introduction}
The study of a pseudo Riemannian manifold $M$ via its isometry group $Iso(M)$ is a central problem in pseudo Riemannian geometry. The larger $Iso(M)$ is, the simpler $M$ is. Many manifolds
have isometry group large enough so that Lie theory can be applied. In mathematics and physics, $n$-dimensional anti de Sitter space is a maximally symmetric Lorentzian manifold with constant negative scalar curvature. In nontransitive cases, $Iso(M)$ is a geometric invariant of $M$
ranking in importance with its curvature and geodesics. This is one
of the reasons that nontransitive actions are of so much interest
to mathematicians. A \co pseudo Riemannian manifold $M$ is an $M$ such that a closed Lie subgroup $G$ of $Iso(M)$ acts on $M$ with a codimension one orbit. Cohomogeneity one Riemannian manifolds have been studied by many mathematicians (see for example, \cite {AA,Ber,PT,PS,S,V}). The problem is still an active one. When the metric is indefinite, there are not so much papers in the literature (see for example \cite{AK1, H^3_1}). In fact there are substantial differences between these two cases. A main difference is that in the Riemannian case, where $G$ is closed in $Iso(M)$, the action is proper, (see \cite {Alee}), which is vital in the study of the subject, while in the indefinite case, this assumption in general does not imply that the action is proper, so the study becomes much more difficult. Also, some of the results and techniques of the definite metric fails for the indefinite metric.

In this paper, which is a continuation of \cite{H^3_1} we study \co three dimensional anti de Sitter space ${\bf adS_3}$. In \cite{H^3_1}, we have studied \co proper actions on ${\bf adS_3}$ and we got some results about the acting group, the orbit space and the causal character of the orbits. Here, the main key of the study is classifying, up to conjugacy, the closed and connected Lie subgroups of $Iso({\bf adS_3})$ which act effectively, isometrically and by \co on ${\bf adS_3}$, in both proper and nonproper cases. When this is done, identifying causal characters of the orbits and the orbit space is an immediate consequence. When the action is proper, all the principal orbits have the same causal character, the same type and there is at most one compact singular orbit (see Theorem \ref{proper}), but in the nonproper case there may be principal orbits of different causal characters, different types and there may be uncountably many singular orbits (see Theorems \ref{3.4.1}, \ref{3.4.2} and \ref{3.4.3}). These are major differences between these two cases.

\section{Preliminaries}

Let $G$ be a Lie group which acts on a connected smooth manifold $M$. The Lie algebra of $G$ is denoted by $\fg$. For each point $x$ in $M$, $G(x)$ denotes the orbit of $x$, and $Stab_G(x)$ is the stabilizer in $G$ of $x$. In this paper, if $x\in M$, then $\fg_x=\{\frac{d}{dt}(\exp(tX)x)|_{t=0}|X\in \fg\}$, so $\fg_x$ does not denote the Lie algebra of the stabilizer in $G$ of $x$, nor does $G_x$ denote the stabilizer in $G$ of $x$. In fact, in this paper the notation $G_x$ will not be used. A smooth manifold $M$ is called of {\it \co} under an action of a Lie
group $G$ if an orbit has codimension one. An action of a Lie group
$G$ on a smooth manifold $M$ is said to be {\it proper} if the mapping
$\varphi:G\times M\rightarrow M\times M , \ (g,x)\mapsto(g.x,x)$ is
proper. Equivalently, for any sequences $x_n$ in $M$ and $g_n$ in $G$, $g_nx_n\rightarrow y$ and $x_n\rightarrow x$ imply that $g_n$ has a convergent subsequence. The $G$-action on $M$ is {\it nonproper} if it is not proper. Equivalently, there are sequences $g_n$ in $G$ and $x_n$ in $M$ such that $x_n$ and $g_nx_n$ converge in $M$ and $g_n\rightarrow \infty$, i.e. $g_n$ leaves compact subsets. For instance, if $G$ is compact, the action is obviously proper. There is a proper action of a Lie
group $G$ on the manifold $M$, if and only if there is a complete
$G$-invariant Riemannian metric on $M$ (see \cite {Alee}). This theorem makes a link between proper actions and Riemannian $G$-manifolds. The orbit
space $M/G$ of a proper action of $G$ on $M$ is Hausdorff, the
orbits are closed submanifolds, and the stabilizers are
compact (see \cite{Adams}). The orbits $G(x)$ and $G(y)$ have the {\it same orbit type} if
$G_x$ and $G_y$ are conjugate in $G$. This defines an equivalence
relation among the orbits of $G$ on $M$. We denote by $[G(x)]$ the
corresponding equivalence class, which is called the {\it orbit type} of
$G(x)$. A submanifold $S$ of $M$ is called a slice at $x$ if there is a $G$-invariant open neighborhood $U$ of $G(x)$ and a smooth equivariant retraction $r:U\rightarrow G(x)$, such that $S=r^{-1}(x)$. A fundamental feature of proper actions is the
existence of slice (see \cite{PT2}), which enables one to define a partial ordering
on the set of orbit types. The partial ordering on the set of orbit types is
defined by, $[G(y)] \leq [G(x)]$ if and only if $G_x$ is conjugate
in $G$ to some subgroup of $G_y$. If $S$ is a slice at $y$, it
implies that $[G(y)] \leq [G(x)]$ for all $x\in S$. Since $M/G$
is connected, there is a largest orbit type in the space of orbit types. Each
representative of this largest orbit type is called a principal
orbit. In other words, an orbit $G(x)$ is principal if and only if
for each point $y\in M$ the stabilizer $Stab_G(x)$ is conjugate to
some subgroup of $Stab_G(y)$ in $G$. Other orbits are called singular.
We say that $x\in M$ is a principal point if $G(x)$ is a principal
orbit.

But for the nonproper action there is not slice in general, so we
can not use the same definitions required the existence of slices as
before, hence we use the definition 2.8.1 of \cite{DK} for
determining the principal, singular or exceptional orbits. According
to it for the action of a Lie group $G$ on the smooth manifold $M$,
The points $x,y \in M$, are said to be of the same type, with
notation $x \thickapprox y$ , if there is a $G$-equivariant
diffeomorphism $\Phi$ from an open $G$-invariant neighborhood $U$ of
$x$ onto an open $G$-invariant neighborhood $V$ of $y$. Clearly this
defines an equivalence relation $\thickapprox$ in $M$. The
equivalence classes will be called {\it orbit types} in $M$, and are
denoted by $M_x^\thickapprox$. If each stabilizer has only
finitely many components, then $x \thickapprox y$ if and only if
$Stab_G(x)$ is conjugate to $Stab_G(y)$ within $G$ and the actions of $Stab_G(x)$, and
$Stab_G(y)$, on $T_xM/T_xG(x)$, and $T_yM/T_yG(y)$, respectively, are
equivalent via a linear intertwining isomorphism (see chapter 2 of
\cite{DK}). The orbit $G(x)$ of $x\in M$ is {\it principal} if its type
$M_x^\thickapprox$ is open in $M$. Any non-principal orbit is called
a {\it singular} orbit. A nonprincipal orbit with the same dimension as a
principal orbit is an {\it exceptional} orbit.

Let $\R^4_2$ denote the $4$-dimensional real
vector space $\R^4$ with the scalar product of signature $(2,2)$
defined by the quadratic form $Q(x)=-(x_1x_4 -x_2x_3)$, where
$x=(x_1,x_2,x_3,x_4)$. The anti de Sitter space ${\bf adS_3}=\{x\in\R^4_2\ |\ Q(x)=-1\}$ is identified with the group $SL(2,\R)$ with biinvariant metric, and the well known classification of one parameter subgroups of $SL(2,\R)$ implies a classification of connected subgroups of the isometry group $(SL(2,\R)\times SL(2,\R))/\Z_2$. Each one parameter subgroup of $SL(2,\R)$ is
    conjugate to one of the groups $A=\{\exp(tX)\ |\ t\in\R\}$, $N=\{\exp(tY)\ |\
    t\in\R\}$ or $K=\{\exp(tZ)\ |\ t\in\R\}$, where
    $X=E_{11}-E_{22}$, $Y=E_{12}$ and $Z=E_{21}-E_{12}$ (see \cite[p.436]{He}). The set  $\{X,Y,Z\}$ is a
    basis for $\mathfrak{sl(\rm 2,\R)}$ and we fix this basis throughout the
    paper. And each two dimensional connected closed Lie subgroup of $SL(2,\R)$ is conjugate to
    $AN$ which is isomorphic to $\Af$, the connected component of the group of affine transformations of the real
    line. The notations $\Af:=AN$,
     $A_t:=e^t E_{11}+e^{-t}E_{22}$, $N_t:=I+tE_{12}$, $K_t:=(\cos t) (E_{11}+E_{22})+(\sin t)(E_{21}-E_{12})$ and $F_{t,s}:=A_t+sE_{12}$, where $s,t\in\R$,
     are used throughout the paper.

%

Let $p$ be an element of ${\bf adS_3}$. The point $p$ is called elliptic, parabolic or hyperbolic if $|tr(p)|$ is less than, equal to or greater that $2$, respectively. This classifies the elements into subsets, not subgroups, since these sets are not closed under multiplication. However, if $p\notin \{\pm I\}$, then it is  elliptic, parabolic or  hyperbolic if $p$ is conjugate into an element of $K$, $N\cup(-N)$ or $A\cup(-A)$, respectively. Topologically, since trace is a continuous map, the set of elliptic elements is an open set, as is the set of hyperbolic elements, while the set of parabolic elements is a closed set.
%
%

    For a Lie group $G$, if $H$ is a subgroup of $G$, we
    use the notation ${\rm diag} (H\times H)$ for the subgroup $\{(h,h)\ |\ h\in
    H\} \subset G\times G$. Let $G\subset SL(2,\R)\times SL(2,\R)$ and $\fg \subseteq
\mathfrak{sl}(2,\R)\oplus
    \mathfrak{sl}(2,\R)$. We denote by
    $p_i:G\rightarrow SL(2,\R)$  and $P_i:\fg \rightarrow \mathfrak{sl}(2,\R)$, where $i=1,2$, the projections
    on the first and second factor, respectively.

\section{Lie Groups acting by \co on ${\bf adS_3}$}

    This section is devoted to classify the Lie subgroups of $SL(2,\R)\times SL(2,\R)$, upto conjugacy, acting effectively, isometrically and by
    \co on ${\bf adS_3}$. This classification is used to determine the
    causal character of the orbits and the orbit space in the next sections. Let
    $\iota:G\rightarrow SL(2,\R)\times SL(2,\R)$ be defined by
    $\iota (g_1,g_2)=(g_2,g_1)$, and $H=\iota(G)$. Then the orbit $G(p)$ is isometric to the orbit $H(p)$, for each
    point $p\in {\bf adS_3}$, and the orbit space ${\bf adS_3}/G$ is
    homeomorphic to ${\bf adS_3}/H$. Thus we consider only the Lie
    subgroups of $SL(2,\R)\times SL(2,\R)$ with $\dim
    p_1(G)\geqslant \dim p_2 (G)$ in the following theorem. First we recall the following lemma from \cite{K} which will be used in the proof of Threorem \ref{1}.

\begin{lemma}\label{K}
     Let $G$ be a Lie group acts on a manifold $M$. Let $H$ be a closed Lie subgroup of $G$
     such that $G/H$
     is compact. Then $G$ acts properly
     if and only if $H$ does.
    \end{lemma}

\begin{theorem}\label{1}
    Let ${\bf adS_3}$ be of cohomogeneity one
    under the action of a connected and closed Lie subgroup $G\subset
    Iso(M)$, then

    (i) $G$ is isomorphic to one of the following Lie groups.
    $$\mathbb{T}^2,\quad R^2 ,\quad \Af  ,\quad \Af\times R ,\quad \Af \times \Af ,\quad SL(2,\R).$$

    (ii) the action is proper if and only if $G$ is conjugate to one of the
    following Lie groups within $SL(2,\R)\times SL(2,\R)$.
    $$ A\times K , \quad N\times K
    ,\quad K\times K  , \quad \Af \times \{I\}\quad
      or \quad G_{FK}=\{(F_{t,s},K_t )\ |\ s,t \in \R\}.$$

      (iii) the action is not proper if and only if $G$ is conjugate to one of the
following Lie groups within $SL(2,\R)\times SL(2,\R)$.
    $$ A\times A, \quad N\times N , \quad A\times N ,\quad \Af\times A,\quad \Af\times N, \quad \Af \times \Af ,$$
    $${\rm diag}(Aff_\circ(\R)\times \Af) , \quad {\rm diag} (SL(2,\R)\times SL(2,\R)),$$
     $$G_{FN}=\{(F_{t,s},N_t ) | s,t \in \R\}, \quad
     G_{FF}=\{(F_{t,s},F_{t,s'})|t,s,s'\in\R\}\quad or$$
     $$ G_{FA}=\{(F_{t,s},A_t ) | s,t \in \R\}.$$

\end{theorem}

\vspace{3mm}

 {\bf Proof :} We break the proof into the
    consideration of several different cases as Lemmas \ref{1.1} to \ref{1.n}, and in each lemma we assume that $G$
    is a closed Lie subgroup of $SL(2,\R)\times SL(2,\R)$
    which acts by \co on ${\bf adS_3}$. Theorem \ref{1} is a direct
    consequence of these lemmas.

\begin{lemma}\label{1.1} If $\dim p_1(G)=1$ then $G$ is conjugate to one of the following Lie groups.
$$(i) \ A\times A, \ (ii)\  N\times N ,\ (iii)\  A\times N ,\ (iv)\ A\times K , \ (v)\  N\times K ,\ (vi)\
K\times K .$$ Furthermore, the action is proper iff $G$ is conjugate
to one of the cases $(iv)$ to $(vi)$.
\end{lemma}
{\bf Proof of Lemma \ref{1.1}.} By the assumption, $\dim
p_1(G)\geqslant \dim p_2(G)$ and the action of $G$ on ${\bf adS_3}$ is by
cohomogeneity one. So $\dim p_2(G)=1$. Then there are $V,W \in
    \mathfrak{sl(\rm 2,\R)}$, such that $p_1(G)=\{\exp (tV)\ |\ t\in\R\}$ and $p_2(G)=\{\exp (tW)\ |\ t\in\R\}$. Let $g_t:=\exp (tV)$ and $h_t:=\exp(tW)$,
    for each $t\in\R$. Let $\fg$ denote the Lie algebra of $G$. Then $P_1(\fg)$ and $P_2(\fg)$ are generated by $V$ and $W$, respectively. Since the action of $G$
    should have a two dimensional orbit, so $\fg=P_1(\fg)\oplus P_2(\fg)$, i.e. $\fg=\{(sV,tW)\ |\ s,t\in\R\}$.
    Hence $G=\{(g_s,h_t)\mid s, t \in \R\}$, and so it is conjugate to one of the
groups stated as the cases $(i)$ to $(vi)$ in Lemma \ref{1.1}. If $p_1(G)$ or
$p_2(G)$ is compact, then the action is proper by Lemma \ref{K}. Let $p_1(G)$ and $p_2(G)$ be noncompact. If $V,W$ are conjugate in
    $SL(2,\R)$, the action is not proper. In fact if $V=pWp^{-1}$ for some $p\in M$,
    then $g_t=ph_tp^{-1}$ for each $t\in \R$. Hence the
    stabilizer of $p$ is
    \begin{eqnarray*}Stab_G(p)&=&\{(g_s,h_t)\in G \mid (g_s,h_t).p=p\} \\ &=&\{(g_t,h_t)\in G\mid g_t=ph_tp^{-1}\},\end{eqnarray*}
    which is a noncompact subgroup, so the action is not proper. This shows that, if $G$ is conjugate to $A\times A$ or $N\times N$, then the action is
    not proper. Thus, to complete the proof, we need only to show that the action of $A\times N$ (case $(iii)$) is not proper. Let $(g_n)=((A_n,N_{e^n}))$ and
    $(p_n)=(e^{-n}E_{11}+E_{12}-E_{21})$. Then $p_n\rightarrow E_{12}-E_{21}$ and $g_np_n\rightarrow I$, but $(g_n)$ has no convergent subsequence. (By a
    simple computation it is seen that the action of $A\times N$ is
    free). {\it End of the proof of Lemma \ref{1.1}}.

\begin{lemma}\label{lem1.2} If $\dim p_1(G)=2$ and $\dim(G)= 2$, then $G$ is conjugate to one of the following
Lie groups that have been stated in Theorem \ref{1}.
$$(i)\ G_{FK} ,\quad (ii)\ G_{FN}  , \quad (iii)\ G_{FA} ,\quad (iv)\ \Af \times \{I\}
\ or\quad (v)\ {\rm diag}(\Af\times \Af).$$
 Furthermore, the action is proper iff $G$ is conjugate to either $G_{FK}$ or $\Af\times \{I\}$. \end{lemma}
 {\bf Proof of
Lemma \ref{lem1.2} :} By the fact that $\dim p_1(G)=\dim(G)=2$, we
have $p_1(G)=Aff_\circ(\R)$, up to conjugacy, and the kernel of the
homomorphism $p_1:G\rightarrow p_1(G)$ is discrete. Hence $p_1$ is a
covering map. The Lie group $p_1(G)$ is simply connected, so $p_1$
is one to one. Hence $p_2\circ p_1^{-1}:p_1(G)\rightarrow p_2(G)$ is
a surjective homomorphism. Let $f=d_I(p_2\circ p_1^{-1})$. If $f=0$,
then $G=\Af \times \{I\}$. If $f\neq 0$ then by the fact that
$\{X,Y\}$ is a basis for $p_1(\fg)$ and $\dim p_2(\fg)=1\ or\ 2$,
the map $f$ sends each $tX+sY$ to an element which is in one of the
following forms, depending on $p_2(\fg)$ to be $K$,
$N$, $A$ or $\Af$ respectively,

$(i) \ (a_1t+b_1s)Z$,

$(ii) \ (a_2t+b_2s)Y$,

$(iii)\ (a_3t+b_3s)X$,

$(iv)\ (a_4t+b_4s)X+(ct+ds)Y$,\\
for some fixed real numbers $a_i,\ b_i,\ c\ and \ d$, where
$1\leqslant i\leqslant4$. In each case the relation
$[f(p_1(\fg)),f(p_1(\fg))]\subseteq f(p_1(\fg))$ implies that
$b_i=0$ for each $1\leqslant i\leqslant4$ and $a_4d=1$. Without less
of generality we may assume that $a_i=1$ for each $1\leqslant
i\leqslant4$ (since we need only the image of $f$) and so we get
that $G$ is conjugate to one of the Lie groups $G_{FK}$, $G_{FN}$ ,
$G_{FA}$ or ${\rm diag}(\Af\times \Af)$, respectively. In the last case,
if $c\neq0$ then $p_2(G)$ is conjugate to $\Af$ which shows that
automorphisms of $\Af$ are conjugacies.

Now in each case we investigate that the action is proper or not.

For the case $G=G_{FK}$, we claim that the action is proper. Let
$((F_{{x_n},y_n},K_{x_n}))$ and $(p_n)$ be sequences in $G_{FK}$ and
${\bf adS_3}$ respectively, such that $F_{{x_n},y_n}p_nK_{-x_n}\rightarrow q$
and $p_n\rightarrow p$ for some $p,q\in {\bf adS_3}$. Let $p_n=(p^n_{ij})$,
$p=(p_{ij})$ and $q=(q_{ij})$, where $i=1,2$ and $j=1,2$. Then
\begin{equation}\label{G_{FK1}}
e^{-x_n}K_{-x_n}\left[\begin{array}{c}p^n_{21}\\p^n_{22}\end{array}\right]
\longrightarrow
\left[\begin{array}{c}q_{21}\\q_{22}\end{array}\right],
\end{equation}
and
\begin{equation}\label{G_{FK2}}
e^{x_n}K_{-x_n}\left[\begin{array}{c}p^n_{11}\\p^n_{12}\end{array}\right]+y_nK_{-x_n}\left[\begin{array}{c}p^n_{21}\\p^n_{22}\end{array}\right]
\longrightarrow
\left[\begin{array}{c}q_{11}\\q_{12}\end{array}\right].
\end{equation}
Since $K_{-x_n}$ is a rotation, $e^{-x_n}$ is convergent to some
point by (\ref{G_{FK1}}), and since $q_{21}q_{22}\neq0$, (note that
$\det (q)=1$) this point is nonzero. Hence $(x_n)$ is convergent.
Using this result and (\ref{G_{FK2}}) shows that $(y_n)$ is
convergent. This completes the proof of our claim about the
properness of the action of $G_{FK}$.

If $G=\Af\times \{I\}$, then the action of $G$ reduce to the left
action of $\Af$ on $SL(2,\R)$ which is free and proper obviously.

If $G=G_{FN}$. Let $(g_n)=((F_{n,0},N_{e^n}))$ and
    $(p_n)=(e^{-n}E_{11}+E_{12}-E_{21})$. Then $p_n\rightarrow E_{12}-E_{21}$ and $g_np_n\rightarrow I$, but $(g_n)$ has no convergent subsequence.
     This shows that the action of $G_{FN}$ is not proper on ${\bf adS_3}$.

If $G=G_{FA}$ then $Stab_G(I)=A$, hence the action is not proper. If
$G={\rm diag}(\Af\times \Af)$, then $Stab_G(I)=G$, and so the action is
not proper. {\it End of the proof of Lemma \ref{lem1.2}}.

\begin{lemma}\label{lem1.3} If $\dim p_1(G)= 2$ and $\dim G\geqslant 3$, then $G$ is conjugate
to one of the following Lie groups within $SL(2,\R)$, and the action
of $G$ on $M$ is not proper.
$$\Af\times A,\quad \Af\times N, \quad G_{FF},\quad or\quad \Af\times \Af.$$  \end{lemma}
{\bf Proof of Lemma \ref{lem1.3}.} By the assumption $\dim
p_1(G)\geqslant \dim p_2(G)$. So there are two following cases

{\bf Case 1. $\dim G=3$.} Since $\dim p_1(G)=2$, we may assume that
$p_1(G)=\Af$. Then the kernel of the linear map $P_1:\fg \rightarrow
\mathfrak{aff}(\R)$ is a one dimensional ideal of $\{0\}\oplus
P_2(\fg)$.

{\it Claim 1.} If $\dim P_2(\fg)=1$, then $G$ is conjugate to either
$\Af\times A$ or $\Af\times N$.\\ Since $\dim (\ker P_1)=1$, we have
$\ker P_1=P_2(\fg)$. Hence $\fg=P_1(\fg)\oplus P_2(\fg)$. By the
well known fact about one and two dimensional subgroups of
$SL(2,\R)$, the Lie group $G$ may be conjugate to one of the groups
$\Af\times A$, $\Af\times N$ or $\Af\times K$. But the action of
$\Af\times K$  on ${\bf adS_3}$ is free and so is not by \co, hence the case
$\Af\times K$ can not occur. {\it End of Claim 1.}

{\it Claim 2.} If $\dim P_2(\fg)=2$, then $G$ is conjugate to
$G_{FF}$.\\
By the assumption of this claim $P_2(\fg)=\mathfrak{aff}(\R)$, up to
conjugacy. Since $\ker P_1$ is a one dimensional ideal of
$\{0\}\oplus P_2(\fg)$, so
\begin{equation}\label{kerP1}
\ker P_1=\{0\}\oplus \{tY|t\in \R\}.\end{equation} On the other
hand, $\fg$ is a three dimensional subalgebra of
$$\mathfrak{aff}(\R)\oplus
\mathfrak{aff}(\R)=\{(tX+sY,t'X+s'Y)|t,t',s,s'\in\R\}.$$ Combining
this with relation (\ref{kerP1}) implies that under the projection
$(tX+sY,t'X+s'Y)\mapsto (t,s,s')$,  $\fg$ maps linear isomorphically
onto $\R^3$. Thus $t'=t'(t,s)$ is a linear function
$t':\R^2\rightarrow \R$. So there are fixed real numbers $a$ and $b$
such that $t'(t,s)=at+bs$. Closeness under the bracket of $\fg$
shows that $b=0$. Therefore $\fg$ has the form
$$\mathfrak{h}_a=\{(tX+sY,atX+s'Y)|t,s,s'\in\R\},$$
where $a$ is nonzero. All nonzero $a$ give conjugate
$\mathfrak{h}_a$. Thus $\fg$ is conjugate to $\mathfrak{h}_1$. {\it
End of Claim 2.} {\it End of Case 1.}

{\bf Case 2.} {\bf $\dim G=4$.} Hence
$P_1(\fg)=P_2(\fg)=\mathfrak{aff}(\R)$ and $\ker P_1$ is a two
dimensional ideal of $\{0\}\oplus \mathfrak{aff}(\R)$. So $\ker
P_1=\{0\}\oplus \mathfrak{aff}(\R)$. Thus
$\fg=\mathfrak{aff}(\R)\oplus \mathfrak{aff}(\R)$ which implies that
$G=\Af\times\Af$. {\it End of Case 2.}

In each case, the stabilizer $G_I$ is not compact, hence the
action is not proper. {\it End of the proof of Lemma \ref{lem1.3}}.

\begin{lemma}\label{1.n}
If $\dim p_1(G)=3$ then $G$ is conjugate to ${\rm diag}(SL(2,\R)\times
SL(2,\R))$ and the action of $G$ on ${\bf adS_3}$ is not proper.
\end{lemma}
{\bf Proof of Lemma \ref{1.n}.} Since $\dim P_1(\fg)=3$ so
$P_1(\fg)=\mathfrak{sl({\rm 2,R)}}$. If $\fg =P_1(\fg)\oplus
P_2(\fg)$ then the action of $G$ on ${\bf adS_3}$ will be transitive,
which is in contrast to the \co assumption. Hence there is a Lie
algebra homomorphism $\varphi:P_1(\fg)\rightarrow P_2(\fg)$.
Clearly, $\varphi(P_1(\fg))$ can not be one or two dimensional,
since $P_1(\fg)$ has no nontrivial ideal as the kernel of $\varphi$.
So $\dim\varphi(P_1(\fg))=3$ and hence
$$P_2(\fg)=\varphi(P_1(\fg))=\mathfrak{sl({\rm 2,R)}}$$
On the other hand, $SL(2,\R)$ is a connected semi-simple Lie group,
so
$$Int(\mathfrak{sl({\rm 2,R)}})=Ad(SL(2,\R))=Aut(\mathfrak{sl({\rm
2,R)}})$$ where $Ad(g)$ is the differential at the identity of the
inner automorphism $x\mapsto gxg^{-1}$, for each $g\in SL(2,\R)$
(see \cite[pp.100-102]{Kn}). Hence there exists $p\in SL(2,\R)$ such
that
$$\varphi (X)=pXp^{-1}\quad\quad,\quad\quad \forall X\in \mathfrak{sl({\rm
2,R)}}.$$ So
$$ \fg=(I,p)[{\rm diag}(\mathfrak{sl({\rm 2,R)}})\times \mathfrak{sl({\rm 2,R)}}](I,p^{-1}).$$
Thus $G$ is conjugate to ${\rm diag}(SL(2,\R)\times SL(2,\R))$. {\it End
of the proof of Lemma \ref{1.n}}.\\
{\it End of the proof of Theorem \ref{1}.}

As a consequence of Theorem \ref{1}, one gets the following
corollary.
\begin{cor}
Let ${\bf adS_3}$ be of \co under the isometric action of a closed and
connected subgroup $G\subseteq Iso_\circ({\bf adS_3})$. If the action is
proper, then $G$ is isomorphic to either $R\times \mathbb{T}$ ,
$\Af$ or $\mathbb{T}^2$.
\end{cor}
\section{Causal characters of the orbits}
Assume that the connected and closed Lie subgroup $G$ of
$Iso({\bf adS_3})$ acts isometrically and by \co on ${\bf adS_3}$, we determine
causal characters of the orbits.

The orbit $G(p)$ is said to be Lorentzian, degenerate or space-like
if the induced metric on $G(p)$ is Lorentzian, degenerate or
Riemannian, respectively. It is called time-like or light-like if
each nonzero tangent vector in $T_pG(p)$ is time-like or null,
respectively. The category into which a given orbit falls is called its {\it causal character}.
\subsection{The action is proper}
Let a Lie group $G$ act by \co and properly on a smooth manifold
$M$.  A result by Mostert (see \cite {Mos}), for the compact Lie
groups, and Berard Bergery (see \cite {Ber}), for the general case,
says that the orbit space $M/G$ is homeomorphic to one of the spaces
$$\R \quad, \quad S^1 \quad, \quad [0,+\infty) \quad , \quad
[0,1].$$ In the following theorem we show that the case $[0,1]$ can
not occur, when $M={\bf adS_3}$. Furthermore, we show that the
causal characters of the principal orbits are the same.

\begin{theorem}\label{proper}
Let ${\bf adS_3}$ be of cohomogeneity one under the proper action of a
connected and closed Lie subgroup $G\subset Iso({\bf adS_3})$. Then one of
the following statements holds.

(i) Each orbit is a Lorentzian surface isometric to $\R\times B$,
where $B$ is anti-isometric to $S^1$. The orbit space ${\bf adS_3}/G$ is diffeomorphic to $\R$.

(ii) Each orbit is a Lorentzian surface diffeomorphic to $\R^2$.
The orbit space ${\bf adS_3}/G$ is diffeomorphic to $S^1$.

(iii) Each orbit is a degenerate surface diffeomorphic to $\R^2$. The orbit space ${\bf adS_3}/G$ is diffeomorphic to $S^1$.

(iv) There is a unique singular orbit anti-isometric to $S^1$
(hence it is time-like), and each principal orbit is a Lorentzian
surface isometric to the Lorentzian torus. The orbit space ${\bf adS_3}/G$ is homeomorphic to $[0,+\infty)$.
%

\end{theorem}

{\bf Proof :} By Theorem \ref{1} one gets all closed and connected
Lie subgroups of $Iso({\bf adS_3})$, up to conjugacy, which act properly
and by \co on ${\bf adS_3}$. We break the proof into two cases as Lemmas
\ref{2.1} and \ref{2.2} , and in each lemma we assume that $G$ is a
closed and connected Lie subgroup of $Iso_\circ({\bf adS_3})$  which acts
properly and by \co on ${\bf adS_3}$, clearly $\dim G = 2$,
and $p=(p_{ij})$ is an arbitrary fixed point of ${\bf adS_3}$, where $1\leqslant i,j \leqslant 2$.

\begin{lemma}\label{2.1}
If $\dim p_1(G)=1$, then one of the statements (i) or
(iv) of Theorem \ref{proper} holds. Furthermore, there is a
singular orbit iff $G$ is conjugate to $K\times K$.
\end{lemma}

{\bf Proof of Lemma \ref{2.1} :} By Theorem \ref{1} $G$ is conjugate
to one of the following Lie groups.
$$ A\times K, \quad N \times K, \quad
    K\times K. $$

So to prove the lemma we need only to study the orbits of the action of these three groups.

{\it Case 1.} ${\bf G=A\times K}$. In this case we claim that the assertion $(i)$ of Theorem \ref{proper} satisfies. Let
    $$\Psi_p(t)=\exp(tX)\ p\ \exp (\alpha tZ)$$
where $\alpha$ is an arbitrary fixed real number. Then
\begin{equation}-\det(\frac{d}{dt}\Psi_p(t)|_{t=0})=-\det(X p-\alpha pZ)=
   -(\alpha^2 +2(p_{11}p_{21}+p_{12}p_{22})\alpha -1).\end{equation}

   The quadratic equation $x^2+2(p_{11}p_{21}+p_{12}p_{22})x-1=0$ has two roots, so $-\det(X p-\alpha
   pZ)$ can be negative, zero or positive, for various $\alpha$. This shows that each orbit is a
    Lorentzian surface. The action of $A\times K$ on ${\bf adS_3}$ is free, so $G(p)$ is diffeomorphic to $\R\times S^1$, for each $p\in {\bf adS_3}$. Since $-\det (\frac{d}{dt}|_{t=0}\exp(tX)p)>0$ and $-\det (\frac{d}{dt}|_{t=0}\exp(tZ)p)<0$, the orbit $G(p)$ is isometric to $\R\times B$, where $B$ is anti-isometric to $S^1$. {\it End of Case 1.}

{\it Case 2.} ${\bf G=N \times K}$. A similar discussion to that of the
case 1 shows that in this case the statement $(i)$ of Theorem \ref{proper} occurs as well. {\it End of Case 2.}

{\it Case 3.} ${\bf G=K\times K}$. In this case we claim that the statement $(iv)$ of Theorem \ref{proper} occurs. We have

$G_I=\{\exp(tZ),\exp(tZ)| \ t \in \R \} \cong SO(2),$\\
where $I$ is the identity matrix. So $G(I)$ is a singular orbit
diffeomorphic to $S^1$. Let
$$\Phi_p(t)=\exp(tZ)\ p\ \exp (\alpha tZ).$$
     Then
    $$-\det(\frac{d}{dt}\Phi_p(t)|_{t=0})=-\det(Z p-\alpha pZ)=
    -(\alpha^2 +(\sum_{i,j} p_{ij}^2)\alpha +1).$$
    Since $p\in SL(2,\R)$, $\sum_{i,j} p_{ij}^2 \geqslant 2$. If $\sum_{i,j} p_{ij}^2=2$
    then $G(p)$ is a time-like singular orbit anti isometric to
    $S^1$. By Theorem 3.1 of \cite [p.38]{Br} the singular orbit $G(p)$ is unique, so all such points,
    $p$ where $\sum_{i,j} p_{ij}^2=2$, belong to the orbit $G(I)$. If $\sum_{i,j} p_{ij}^2 > 2$ then
    $-\det(\frac{d}{dt}\Phi_t(p)|_{t=0})$
    can be positive, zero or negative for different $\alpha$, so $G(p)$ is a Lorentzian
    principal orbit. The Lie group $G$ is isomorphic to $\mathbb{T}^2$, so $G(p)$ is isometric to the Lorentzian torus.{\it End of Case 3.}\\
    {\it End of the proof of Lemma \ref{2.1}.}
\begin{lemma}\label{2.2}
If $\dim p_1(G)\geqslant 2$, then one of the statements (ii) or (iii)
of Theorem \ref{proper} holds. Furthermore there are some light-like
orbits if and only if $G$ is conjugate to $G_{FK}$.
\end{lemma}

{\bf Proof of Lemma \ref{2.2} :} By Theorem \ref{1}, $G$ is
conjugate to one of the Lie groups $\Af \times \{I\}$ or  $G_{FK}$.
To prove the lemma we need only to study the actions of these groups on ${\bf adS_3}$.

{\it Case 1.} ${\bf G=G_{FK}}$. We claim that, this case leads to the assertion {\it (ii)} of Theorem \ref{proper}. Since $G_{FK}$ is isomorphic to $\Af$ and the action is proper, so $Stab_G(p)=\{I\}$. This shows that the action is free. The set $\{(X,Z) , (Y,0) \}$ is a basis for the Lie algebra $\fg$. To find the causal
   character of the orbit $G(p)$, let
   $$\Phi_p(t)=\exp (t(X+\alpha Y)) p \exp (-tZ)$$
   where $\alpha$ is an arbitrary fixed real number. Then
   \begin{equation*}-\det(\frac{d\Phi_p}{dt}(0))=-\alpha(p_{21}^2+p_{22}^2)+2(p_{11}-p_{12})p_{22}.\end{equation*}
   Hence $-\det(\frac{d\Phi_p}{dt}(0))$ may be negative, zero or positive for various $\alpha$, so
   $G(p)$ is a Lorentzian surface diffeomorphic to $\R^2$. {\it End of Case 1.}

{\it Case 2.} ${\bf G=\Af\times \{I\}}$. In this case we claim that the statement $(iii)$ of Theorem \ref{proper} holds. By the fact that $\{ X , Y \}$ is a
basis for the Lie algebra $\fg$, if one defines
    $\Psi_p(t)=\exp(tX)p$ and $\Phi_p(t)=\exp(tY)p$.
   Then the relations

   $\frac{d\Psi_p}{dt}(t)=Xp \lo -\det(Xp)=1,$

   $\frac{d\Phi_p}{dt}(t)=Yp \rightarrow -\det(Yp)=0,$\\
   show that the orbit $G(p)$ is a degenerate orbit. Since, in this case, the action of $G$ on ${\bf adS_3}$ is free, the orbit $G(p)$ is diffeomorphic to $\R^2$. {\it End of Case 2.}\\
   {\it End of the proof of Lemma \ref{2.2}}.\\
   {\it End of the proof of Theorem \ref{proper}}.

As a consequence of Theorem \ref{proper} and its proof we have the
following corollary.
\begin{cor}
Let ${\bf adS_3}$ be of cohomogeneity one under the proper action of a
connected and closed Lie subgroup $G\subset Iso({\bf adS_3})$. Then the
following assertions hold.

(a) There is no space-like orbit.

(b) There is a degenerate orbit if and only if each orbit is a
degenerate surface diffeomorphic to $\R^2$, if and only if $G$ is
conjugate to $\Af\times \{I\}$.
\end{cor}

\subsection{The action is not proper}
Let $G$ be a closed and connected Lie subgroup of $Iso({\bf adS_3})$ which acts nonproperly and by \co on ${\bf adS_3}$. In this subsection we determine the causal character of the orbits. Throughout of the section it is assumed that $p=(p_{ij})$ is an arbitrary fixed point of ${\bf adS_3}$, where $1\leqslant i,j \leqslant 2$.
\begin{theorem}\label{3.4.1}
Let ${\bf adS_3}$ be of cohomogeneity one under the action of a connected
and closed Lie subgroup $G\subset Iso({\bf adS_3})$. If the action is not
proper and $\dim p_1(G))=1$, then $G$ is conjugate to one of the Lie groups $ A\times A$, $A\times N$ or $N\times N$. Furthermore, the following assertions hold.

(i)  If $G$ is conjugate to $A\times A$, then there are four
space-like singular orbits diffeomorphic to $\R$ of the same type.
Each principal orbit is space-like or Lorentzian and all of them are
of the same type. The union Lorentzian principal orbits is open in
${\bf adS_3}$.

(ii) If $G$ is conjugate to $A\times N$, then the action is free and each orbit is principal. All of them are of the same orbit type. There are just four degenerate orbits and each other orbit is Lorentzian.

(iii) If $G$ is conjugate to $N\times N$, then there are uncountably many light-like singular orbits of the same type, where each of them is diffeomorphic to $\R$. Each principal orbit is a Lorentzian surface diffeomorphic to $\R^2$, and all of them are of the same type.
\end{theorem}
{\bf Proof of Theorem \ref{3.4.1} :} By Theorem \ref{1} $G$ is
conjugate to one of the following Lie groups.
$$ A\times A ,\quad A\times N, \quad N\times N.  $$

{\bf (i)}  ${\bf G=A\times A}$. Then $G(p)=\{\exp (tX) p \exp (-sX)
| s,t\in \R\}.$ By a simple computation it is seen that $Stab_G(p)$
is conjugate to ${\rm diag}(A\times A)$ if  $p\in \pm(A\cup JA)$,
and $Stab_G(p)=\{I\}$ if  $p\notin \pm(A\cup JA)$, where
$J=E_{12}-E_{21}$. This shows that the orbit $G(p)$ is singular if
and only if $p\in \pm(A\cup JA)$. On the other hand, $p\in A$
implies that $G(p)=A$, and $p\in JA$ implies that $G(p)=JA$. Hence
there are just four singular orbits corresponding to the points $\pm
I$ and $\pm J$. These singular orbits are of the same type, since
the stabilizer of each point of them is conjugate to ${\rm
diag}(A\times A)$.

To determine the causal character of the orbits, let $$\Psi_p(t)=\exp (tX) p \exp(\alpha tX),$$ where $\alpha$ is an arbitrary fixed real number.
Then
$$-\det(\frac{d\Psi_p}{dt}(0))=\alpha^2-2(2p_{11}p_{22}-1)\alpha+1.$$
Hence the following cases occur.

(a) If $p\in \pm(A\cup JA)$ (and so $p_{11}p_{22}$ is $0$ or $1$),
then $-\det(\frac{d\Psi_p}{dt}(0))=(\alpha\pm1)^2$. Hence the orbit
$G(p)$ is a one dimensional space-like orbit diffeomorphic to $\R$
(note that the directions where $\alpha=\mp1$ belong to $\fg_p$).

(b)  If $p\notin \pm(A\cup JA)$ and $0\leqslant
p_{11}p_{22}\leqslant 1$, then the orbit $G(p)$ is a space-like
surface diffeomorphic to $\R^2$. So the union of space-like
principal orbits, for which each point of them has no zero entry, is
open in ${\bf adS_3}$.

(c) If $p_{11}p_{22}<0$ (or $p_{11}p_{22}>1$), then $G(p)$ is a
Lorentzian surface diffeomorphic to $\R^2$. So the union of
Lorentzian principal orbits is open in ${\bf adS_3}$.

{\bf (ii)} ${\bf G= A\times N}$. It is easily seen that the action of $G$ on ${\bf adS_3}$ is free, and so there is no singular orbit and each orbit is diffeomorphic to $\R^2$. Hence there is only one principal orbit type.
To determine the causal character of the orbits, let
    $$\Psi_p(t)=\exp(tX) p \exp (\alpha tY),$$ where $\alpha$ is an arbitrary fixed real number.
Then
    $$-\det(\frac{d}{dt}\Psi_p(t)|_{t=0})=-\det(X p-\alpha pY)=1+2\alpha p_{11}p_{21}.$$
    Hence $p_{11}p_{21}\neq 0$ implies that the polynomial $1+2\alpha p_{11}p_{21}$ can be
    positive, zero or negative for various $\alpha$, which shows that the orbit $G(p)$
    is a Lorentzian surface. If $p_{11}p_{21}=0$ then $-\det(\frac{d}{dt}\Psi_p(t)|_{t=0})=1$
    for each $\alpha \in \R$. On the other hand, $(0,Y)$ is a null direction tangent to $G(p)$ in $p$, so the orbit $G(p)$ is a degenerate principal
    orbit when $p_{11}p_{21}=0$.

    We claim that there are just four degenerate orbits. To prove our claim, first suppose that $p_{11}=0$. Then $p'\in G(p)$, if and only if there are some real numbers $s$ and $t$
    so that $A_t p N_{-s}=p'$. Solving this equation leads one to get that there are such real numbers $s$ and $t$ if and only if $p_{11}'=0$, $p_{12}p_{12}'>0$.
    A similar discussion shows that for the case $p_{21}=0$, we have $p'\in G(p)$ if and only if $p_{21}'=0$, $p_{11}p_{11}'>0$.
    Since $G(p)$ is a degenerate orbit if and only if
    $p_{11}p_{21}=0$, our claim is proved.

{\bf (iii)} ${\bf G=N\times N}$. Then for each point $p\in {\bf adS_3}$, where $p_{21}=0$, we have $Stab_G(p)=\{(N_t,N_{s(t)})|t\in\R,\ s(t)=\frac{p_{22}}{p_{11}}t\}$, and for each other point we have $Stab_G(p)=\{I\}$. This shows that $G(p)$ is a singular orbit if and only if $p_{21}=0$. There are uncountably many singular orbits, since (for $p$ where $p_{21}=0$)
$$G(p)=\{p_{11}E_{11}+(p_{12}+tp_{22}-sp_{11})E_{12}+p_{22}E_{22}|s,t\in\R\}.$$ Comparing the stabilizers shows that there is one principal orbit type and one singular orbit type, since the stabilizer of each singular point is conjugate to ${\rm diag}(N\times N)$.

To determine the causal character of the orbits, let $\Psi_p(t)=\exp (tY) p \exp(\alpha tY)$, where $\alpha$ is an arbitrary fixed real number. So
$$-\det(\frac{d\Psi_p}{dt}(0))=\alpha p_{21}^2.$$ Hence $G(p)$ is a
Lorentzian orbit if $p_{21}\neq 0$, and it is light-like if $p_{21}=0$. Thus each principal orbit is a Lorentzian surface diffeomorphic to $\R^2$ and each singular orbit is a one dimensional light-like submanifold diffeomorphic to $\R$. {\it End of the proof of Theorem \ref{3.4.1}}.

\begin{theorem}\label{3.4.2}
Let ${\bf adS_3}$ be of cohomogeneity one under the action of a connected
and closed Lie subgroup $G\subset Iso({\bf adS_3})$. Let the action is not
proper. If $\dim p_1(G))=1$ and $\dim G=2$, then $G$ is conjugate to one of the Lie groups ${\rm diag}(\Af \times \Af)$, $G_{FN}$ or $G_{FA}$. Furthermore, the following assertions hold.

(i)  If $G$ is conjugate to ${\rm diag}(\Af \times \Af)$, then there
are two zero dimensional singular orbits and uncountably many one
dimensional light-like singular orbits, where each of them is
diffeomorphic to $\R$. There are three types of singular orbits.
Each principal orbit may be space-like, degenerate or Lorentzian
surface diffeomorphic to $\R^2$, but all of them are of the same
type.

(ii) If $G$ is conjugate to $G_{FN}$, then the action is free and
each orbit is principal. Two of the orbits are degenerate and each
other orbit Lorentzian. All of them are of the same orbit type.

(iii) If $G$ is conjugate to $G_{FA}$, then there are uncountably many light-like singular orbits diffeomorphic to $\R$ of the same type. Each principal orbit is a degenerate or Lorentzian surface diffeomorphic to $\R^2$. All the principal orbits are of the same type.
\end{theorem}
{\bf Proof of Theorem \ref{3.4.2} :}
By Theorem \ref{1} $G$ is
conjugate to one of the following Lie groups.
$${\rm diag}(\Af \times \Af), \quad G_{FN},\quad G_{FA}.$$ 
So we consider the actions of these groups on ${\bf adS_3}$ as follows.

{\bf (i)} ${\bf G={\rm diag}(\Af\times \Af).}$ The set
$\{(X,X),(Y,Y)\}$ is a basis for $\fg$, so to determine the causal
character of the orbits, let
 $$\Phi_p(t)=\exp t(\alpha X+\beta Y)p\exp
(-t(\alpha X+\beta Y)),$$ where $\alpha$ and $\beta$ are fixed real
numbers. Then
\begin{equation}\label{diag(Af)}
-\det(\frac{d\Phi_p}{dt}(0))=\beta^2 p_{21}^2-2\alpha p_{21}(2\alpha
p_{12}+\beta(p_{22}-p_{11})).\end{equation}

Let $p\notin Aff(\R)$, (equivalently $p_{21}\neq 0$). Then
$Stab_G(p)=\{I,I\}$. So by using (\ref{diag(Af)}) and discussing on
different values of $p=(p_{ij})$, one gets that $G(p)$ may be a
space-like, degenerate or Lorentzian surface diffeomorphic to
$\R^2$. Each $G(p)$ is a principal orbit, since $M_p^\thickapprox$
is open in ${\bf AdS_3}$. All of these orbits, $G(p)$'s where
$p\notin Aff(\R)$, are of the same type since their stabilizers are
identity.

Let $p\in Aff(\R)$ (equivalently $p_{21}=0$). Then
\begin{equation}\label{G-diag(Af)} G(p)=\{p_{11}E_{11}+((p_{22}-p_{11})s+p{12}e^t)e^tE_{12}+p_{22}E_{22}|s,t\in\R\},\end{equation}
which is a singular orbit and it is light-like by (\ref{diag(Af)}).
Thus $G(p)$ is a singular orbit if and only if $p\in Aff(\R)$. To
determine the types of orbits we consider three cases as follows.

{\it Case} (a). $p\in Aff(\R)-(N\cup (-N))$. Then
$$Stab_G(p)=\{(F_{t,s(t)},F_{t,s(t)})| t\in \R,\ s(t)=\frac{p_{12}(e^{-t}-e^t)}{p_{22}-p_{11}}\},$$ which is conjugate to ${\rm diag}(A\times
A)$. By using (\ref{G-diag(Af)}), one gets that $p'\in G(p)$ if and
only if $p'_{21}=0$ and $p_{11}=p'_{11}$. This implies that there
are uncountably many one dimensional singular orbits. These orbits
are of the same type, since the stabilizer of any point of them is
conjugate to ${\rm diag}(A\times A)$.

{\it Case} (b). $p\in N\cup(-N)-\{\pm I\}$. By using
(\ref{G-diag(Af)}), we have $p'\in G(p)$ if and only if $p'_{21}=0$,
$p_{11}=p'_{11}$ and $p_{12}p'_{12}>0$. Since $p_{11}=\pm 1$, there
are just four singular orbits in this case. The stabilizer of any
point of these four orbits is ${\rm diag}(N\times N)$. So these
orbits are of the same type, but this type differs from that of in
case (a).

{\it Case} (c). $p\in \{\pm I\}$. Then $p$ is fixed by $G$ and this
case yields two zero dimensional singular orbits with the same type.
Obviously, this type is different from those of in cases (a) and
(b).

%
%

 {\bf (ii)}  ${\bf G=G_{FN}}$. By a simple computation one gets that the action is free, so there is no singular orbit. Since $(Y,0)\in \fg$ is a null direction, so there is no space-like orbit. Any other direction in $\fg$ may be obtained by $(X,Y)+\beta(Y,0),$ where $\beta$ is an arbitrary fixed real number. Let
     $$\Phi_p(t)=\exp t(X+\beta Y) p \exp (-tY).$$  Then
     $$-\det (\frac{d\Phi_p}{dt}(0))= 1+2p_{11}p_{21}+\beta p_{21}^2,$$
     which shows that $G(p)$ is a Lorentzian (resp. degenerate) principal
     orbit if $p_{21}\neq 0$ (resp. $p_{21}=0$). If $p_{21}=0$, Then
     $p'\in G(p)$ if and only if $p'_{21}=0$ and $p_{11}p'_{11}>0$.
     This shows that only two of the orbits are degenerate.

{\bf (iii)}  ${\bf G=G_{FA}}$. If $p\notin Aff(\R)$ then
$Stab_G(p)=\{I\}$, and if $p\in Aff(\R)$ (i.e. $p_{21}= 0$) then
  \begin{equation}\label{G_FA1} Stab_G(p)=\{(F_{t,s(t)},A_t)|t\in\R,\ s(t)=\frac{p_{12}(e^t-e^{-t})}{p_{22}}\},\end{equation}
  which is conjugate to ${\rm diag}(A\times A)$. Hence $G(p)$ is a singular orbit if and only if $p\in Aff(\R)$. If $p_{21}=0$, then
  it is easily seen that $p'\in G(p)$ if and only if $p'_{21}=0$ and $p'_{22}=p_{22}$. So there are uncountably many singular orbits diffeomorphic to $\R$ of the same type.

  The set $\{(X,X),((0,Y)\}$ is a basis for $\fg$. So to find the causal character of the orbits, let
 $$\Phi_p(t)=\exp (t(\alpha X+\beta Y)) p \exp
(-t\alpha X),$$ where $\alpha$ and $\beta$ are
fixed real numbers. Then
\begin{equation*}\label{G_{FA}}
-\det(\frac{d\Phi_p}{dt}(0))=-2\alpha p_{21}(2\alpha p_{12}+\beta p_{22}).\end{equation*}
 Combining this with (\ref{G_FA1}), one gets that $G(p)$ is a one dimensional light-like submanifold, a Lorentzian surface or a degenerate surface, if $p_{21}=0$, $p_{21}p_{22}\neq0$, or $p_{21}\neq 0$ and $p_{22}=0$, respectively. Thus each principal orbit is a degenerate or Lorentzian surface diffeomorphic to $\R^2$, and all of them are of the same type. {\it End of the Proof of Theorem \ref{3.4.2}}.

\begin{theorem}\label{3.4.3}
Let ${\bf adS_3}$ be of cohomogeneity one under the action of a connected
and closed Lie subgroup $G\subset Iso({\bf adS_3})$. Let the action is not
proper. If $\dim G\geqslant 3$, then $G$ is conjugate to one of the following Lie groups
$$(i)\Af\times A,\ (ii)\Af\times N, \ (iii) \Af \times \Af ,$$
    $$ (iv) {\rm diag} (SL(2,\R)\times SL(2,\R)),\ (v) G_{FF} .$$ Furthermore, the following assertions hold.

(i)  If $G$ is conjugate to $(\Af \times A)$, then there are only
four degenerate exceptional orbits diffeomorphic to $\R^2$, and there are four
orbits which are open submanifolds. All exceptional orbits are of the
same type.

(ii) If $G$ is conjugate to $(\Af \times N)$, then there are only
two degenerate exceptional orbits diffeomorphic to $\R^2$, and there are two
orbits which are open submanifolds. All exceptional orbits are of the
same type.

(iii) If $G$ is conjugate to $(\Af \times \Af)$, then there are only
two degenerate exceptional orbits diffeomorphic to $\R^2$ and two
orbits which are open submanifolds. All exceptional orbits are of the
same type.

(iv) If $G$ is conjugate to ${\rm diag}(SL(2,\R)\times SL(2,\R))$,
then the orbit of each elliptic element is a space-like principal
orbit diffeomorphic to $\R^2$, the orbit of each parabolic element
is a space-like principal orbit diffeomorphic to $\R\times S^1$, and
the orbit of each hyperbolic element is a Lorentzian principal orbit
diffeomorphic to $\R\times S^1$. There are three orbit types
corresponding to the orbits of elliptic, parabolic and hyperbolic
elements. There are two zero dimensional singular orbits
corresponding to $\{\pm I\}$.

(v)  If $G$ is conjugate to $G_{FF}$, then there are uncountably
many light-like singular orbits of the same type, and two orbits
which are open submanifolds. Each singular orbit is diffeomorphic to
$\R$.
\end{theorem}
{\bf Proof of Theorem \ref{3.4.3} :} By Theorem \ref{1} $G$ is
conjugate to one of the mentioned Lie groups. Let $(g_{ij}^p)$ denote a typical element of $G(p)$, where $1\leqslant i,j\leqslant 2$.

    {\bf (i)} ${\bf G= \Af \times A}$. The set $\{(X,0),(0,Y),(0,X)\}$ is a basis for $\fg$. To determine the causal character of the orbits, suppose that
    $$\Phi_p(t)=\exp t(X+\alpha Y) p\exp (-\beta tX),$$ where $\alpha$ and $\beta$ are fixed real numbers. Then
    \begin{equation}\label{eq1}-\det(\frac{d\Phi_p}{dt}(0))=\beta^2+2(p_{11}p_{22}+p_{12}p_{21})\beta +1+2\alpha \beta p_{21}p_{22}.\end{equation}

    If $p_{21}=0$ then $Stab_G(p)=\{(F_{t,s(t)},A_{t} )|t \in \R,\ s(t)=-\frac{p_{12}}{p_{22}}(e^t-e^{-t})\}$, and so $G(p)$ is a two dimensional, diffeomorphic to $\R^2$, degenerate submanifold of ${\bf adS_3}$ by (\ref{eq1}).  By a simple computation one sees that $g_{21}^p=0$ and $p_{22}g_{22}^p>0$. This shows that there are only two
    orbits $G(p)$ with $p_{21}=0$. Hence $M_p^\thickapprox$ is not open in ${\bf adS_3}$, and so $G(p)$ is an exceptional orbit.

     If $p_{21}\neq 0$ and $p_{22}=0$ then $Stab_G(p)=\{(F_{t,s(t)},A_{-t} )\ |t \in \R,\ s(t)=-\frac{p_{11}}{p_{21}}(e^t-e^{-t})\}$, and
     so $G(p)$ is diffeomorphic to $\R^2$, which is a degenerate surface by (\ref{eq1}).
     Let $p'\in {\bf adS_3}$. Then $p'\in G(p)$ if and only if there are some real numbers
    $s,t$ and $u$ such that $F_{t,s}pA_{-u}=p'$. By a simple computation, one gets that there are such real numbers $s,t$ and $u$
    if and only if  $p_{21}p_{21}'>0$ and $p_{22}'=0$. Thus there are only two orbits with $p_{21}\neq 0$ and $p_{22}=0$, which shows that $M_p^\thickapprox$ is not open in ${\bf adS_3}$. Hence $G(p)$ is an exceptional orbit. Furthermore,
 these two orbits are of the same type, since $gStab_G(p)g^{-1}=Stab_G({p')}$, where $g=(N_{(p_{11}'/p_{21}')-(p_{11}/p_{21})},I)$.

      If $p_{21} p_{22}\neq0$ then $Stab_G(p)=\{I\}$, and so $\dim G(p)=3$.
    Hence $G(p)$ is an open submanifold of ${\bf adS_3}$. Let $p'\in {\bf adS_3}$. By a similar discussion of that of the preceding case, one gets that $p'\in {\bf adS_3}$ if and only if $p_{21}p_{21}'>0$ and $p_{22}p_{22}'>0$. This shows that there are exactly four disjoint open
    submanifolds as the orbits.

{\bf (ii)} ${\bf G= \Af \times N}$. The set $\{(X,0),(Y,0),(0,Y)\}$
is a basis for $\fg$. Let $$\Phi_p(t)=\exp
    t(X+\alpha Y) p \exp (-\beta tY),$$ where $\alpha$ and $\beta$ are fixed real numbers. Then
    \begin{equation}\label{eq2}-\det(\frac{d\Phi_p}{dt}(0))=\alpha\beta p_{21}^2+2\beta p_{11}p_{21} +1.\end{equation}
     Hence one may consider the following cases.

    (a) If $p\notin Aff(\R)$, i.e. $p_{21}\neq 0$, then the action is free, and so $G(p)$ is an
    open submanifold of ${\bf adS_3}$ diffeomorphic to $\R^3$. Let $p'\in {\bf adS_3}$. Then $F_{t,s}pN_{u}^{-1}=p'$ if and only if
    $$t=\ln \frac{p_{21}}{p_{21}'}\quad ,\quad s=\frac{p_{11}'p_{21}'-p_{11}p_{21}}{p_{21}p_{21}'} \quad
    ,\quad u=\frac{p_{22}p_{21}'-p_{22}'p_{21}}{p_{21}p_{21}'}.$$
    Hence $G(p)=G(p')$ if and
    only if $p_{21}p_{21}'>0$. Thus there are two open submanifolds as the
    orbits.

    (b) If $p\in Aff(\R)$, i.e. $p_{21}=0$, then $Stab_G(p)=\{(N_{s},N_{-\frac{p_{22}}{p_{11}}s}) |
    s\in\R\}$ and so $G(p)$ is a degenerate (not light-like) two dimensional submanifold by
    (\ref{eq2}). In this case, if $p'$ is another point of $ {\bf adS_3}$ with $p_{21}'=0$ and $p_{11}p_{11}'>0$, then with the real
    numbers $t=\ln \frac{p_{11}'}{p_{11}}$ and $s=p_{11}p_{12}'-p_{11}'p_{12}$ one gets that
    $F_{t,s}p=p'$. Hence $G(p)=G(p')$ if and only if $p_{21}=p_{21}'=0$ and
    $p_{11}p_{11}'>0$ (equivalently $p_{22}p_{22}'>0$). Thus there are only two distinct
    orbits with $p_{21}=0$. these two orbits are of the same type,
    since their stabilizers are conjugate.

{\bf (iii)} ${\bf G=\Af \times \Af}$. The set
$\{(X,0),(Y,0),(0,X),(0,Y)\}$ is a basis for $\fg$. Let
$$\Phi_p(t)=\exp t(\alpha X+\beta Y)p\exp -t(\gamma X+\eta Y),$$
where $\alpha,\beta,\gamma$ and $\eta$ are fixed real numbers. Then
\begin{equation}\label{eq3}-\det(\frac{d\Phi_p}{dt}(0))=\beta\eta p_{21}^2+2\alpha\beta
p_{21}p_{22}-2\gamma\eta p_{11}p_{21}-4\alpha\gamma
p_{11}p_{22}+(\alpha+\gamma)^2.\end{equation}

If $p\in Aff(\R)$, i.e. $p_{21}=0$, then
$Stab_G(p)=\{(F_{t,s},F_{t'(t),s'(t,s)})\in G| t'(t)=t, \
s'(t,s)=\frac{p_{12}(e^t-e^{-t})+p_{22}s}{p_{11}}\}$, and
$-\det(\frac{d\Phi_p}{dt}(0))=(\alpha-\gamma)^2$ by (\ref{eq3}).
Hence $G(p)$ is a two dimensional degenerate orbit. Let $p'$ be an
arbitrary point of ${\bf adS_3}$. Then $p'\in G(p)$ if and only if there
are some real numbers $t,s,t'$ and $s'$ such that
$F_{t,s}pF_{t',s'}^{-1}=p'$, and this holds if and only if
$p_{21}'=0$ and $p_{11}p_{11}'>0$ (equivalently $p_{22}p_{22}'>0$).
This shows that there are only two distinct orbits with $p_{21}=0$.
They are of the same orbit type, since their stabilizers are
conjugate to ${\rm diag}(\Af\times \Af)$.

If $p\notin Aff(\R)$, i.e., $p_{21}\neq 0$, then
$$G_p=\{(F_{t,s(t)},F_{t'(t),s'(t)})\in G\ |\ t'(t)=-t\ ,\
s(t)=\frac{p_{11}(e^{-t}-e^{t})}{p_{21}} ,\
s'(t)=\frac{p_{22}(e^{-t}-e^t)}{p_{21}}\}.$$ Hence $\dim G(p)=3$ and
so $G(p)$ is an open submanifold of ${\bf adS_3}$ diffeomorphic to $\R^3$.
If $p'$ is another point of ${\bf adS_3}$, then $p'\in G(p)$ if
$p_{21}p_{21}'>0$. On the other hand, if $p_{21}p_{21}'>0$, then
with the real numbers
$$ t'=0\ ,\quad t=\ln \frac{z}{z'}\ , \quad s=\frac{x'z'-xz}{zz'}\ ,\quad s'=\frac{wz'-w'z}{zz'}\ ,$$
one gets that $F_{t,s}pF_{t',s'}^{-1}=p'$. This shows that there are
only two orbits which are open submanifolds.

{\bf (iv)} ${\bf G={\rm diag}(SL(2,\R)\times SL(2,\R))}$. In this
case $Stab_G(p)=Cent_G(p)$. For each $p'\in {\bf adS_3}$,  we have
$G(p)=G(p')$ if and only if $p$ is conjugate to $p'$ within
$SL(2,\R)$. So to study the causal character of the orbit of each
point, we need only to consider points upto conjugacy. Let
$$\Phi_p(t)=\exp t(\alpha X+\beta Y+\gamma Z) p \exp(-t(\alpha X+\beta Y+\gamma
Z)),$$ where $\alpha,\ \beta$ and $\gamma$ are arbitrary fixed real
numbers. Then
\begin{equation}\label{slsl}-\det(\frac{d\Phi_p}{dt}(0))=(\gamma
p_{12}-(\beta-\gamma)p_{21})^2-(2\alpha
p_{12}+(\beta-\gamma)(p_{22}-p_{11}))(\gamma (p_{22}-p_{11})+2\alpha
p_{21}).
\end{equation}
If $p\in \{I,-I\}$, then $G(p)=\{p\}$.\\ Let $p\notin \{I,-I\}$.

 If $p$ is an elliptic, parabolic or hyperbolic element, then
$Stab_G(p)$ is equal to $K$, $N\cup (-N)$ or $A\cup (-A)$,
respectively.

If $p$ is an elliptic element,then by (\ref{slsl}), on gets that
$-\det(\frac{d\Phi_p}{dt}(0))=(\beta^2+4\alpha^2)p_{21}$,  which
implies that $G(p)$ is a space-like orbit diffeomorphic to $\R^2$.

If $p$ or $-p$ is a parabolic element, then $-\det(\frac{d\Phi_p}{dt}(0))=\gamma^2 p_{12}^2$ by (\ref{slsl}), so $G(p)$ is a space-like orbit diffeomorphic
 to $\R\times S^1$.

If $p$ or $-p$ is a hyperbolic element, then $-\det(\frac{d\Phi_p}{dt}(0))=\gamma(\beta-\gamma)(p_{22}-p_{11})^2$ by (\ref{slsl}), which shows that  $G(p)$ is a
Lorentzian orbit diffeomorphic to $\R \times S^1$.\\
Comparing the stabilizer subgroups shows that there are three orbit
types.

 {\bf (v)} ${\bf G=G_{FF}}$.  Then by a simple computation one gets that
$$\begin{array}{l}
    p_{21}\neq 0 \Rightarrow Stba_G(p)=\{I\}, \\
  p_{21}=0 \Rightarrow Stba_G(p)=\{(F_{t,s},F_{t,s'(t,s)})|t,s\in\R,\ s'(t,s)=\frac{y(e^t-e^{-t})+p_{22}s}{p_{11}}\}.
\end{array}$$
So an orbit $G(p)$, where $p_{21}\neq 0$, is an open submanifold. To
determine causal characters of singular orbits we do as follows. The
set $\{(X,0),(Y,0),(0,Y)\}$ is a basis for $\fg$. Let
 $$\Phi_p(t)=\exp (t(\alpha X+\beta Y)) p \exp
(-t(\alpha X+\gamma Y)),$$ where $\alpha$, $\beta$ and $\gamma$ are
fixed real numbers. Then
\begin{equation}\label{G_{FF}}
-\det(\frac{d\Phi_p}{dt}(0))=-\beta\gamma p_{21}^2+4\alpha^2
p_{12}p_{21}-2\alpha\gamma p_{11}p_{21}+2\alpha\beta
p_{21}p_{22},\end{equation}
which shows that any singular orbit $G(p)$, where $p_{21}=0$, is a one dimensional light-like submanifold.\\
Let $p_{21}\neq 0$. By a similar discussion of that of Case (iii),
one gets that $p'\in G(p)$ if and only if $p_{21}p'_{21}>0$. This
implies that there are just two open submanifolds as the orbits. If
$G(p)$ is a singular orbit, then $G(p)=G(p')$ if and only if
$p_{11}=p'_{11}$. Hence there are uncountably many singular orbits.

Comparing the stabilizers shows that there is one singular orbit
type. {\it End of the proof of Theorem \ref{3.4.3}}.

As a consequence of Theorems \ref{3.4.1}, \ref{3.4.2} and \ref{3.4.3}, one gets the following corollary.

\begin{cor}
Let ${\bf adS_3}$ be of cohomogeneity one under the action of a connected
and closed Lie subgroup $G\subset Iso({\bf adS_3})$. If the action is not
proper, then the following statements hold.

(i) There is no space-like exceptional orbit.

(ii) There is some exceptional orbit if and only if there is some open submanifold as an orbit.

(iii) All exceptional orbits are of the same type.

(iv) If there is an exceptional orbit, then there is no principal orbit.

(v) There are at most four exceptional orbits.
\end{cor}

\section{The orbit space}

\begin{figure}‎
‎\centering‎
‎\includegraphics[width=1.43\columnwidth]{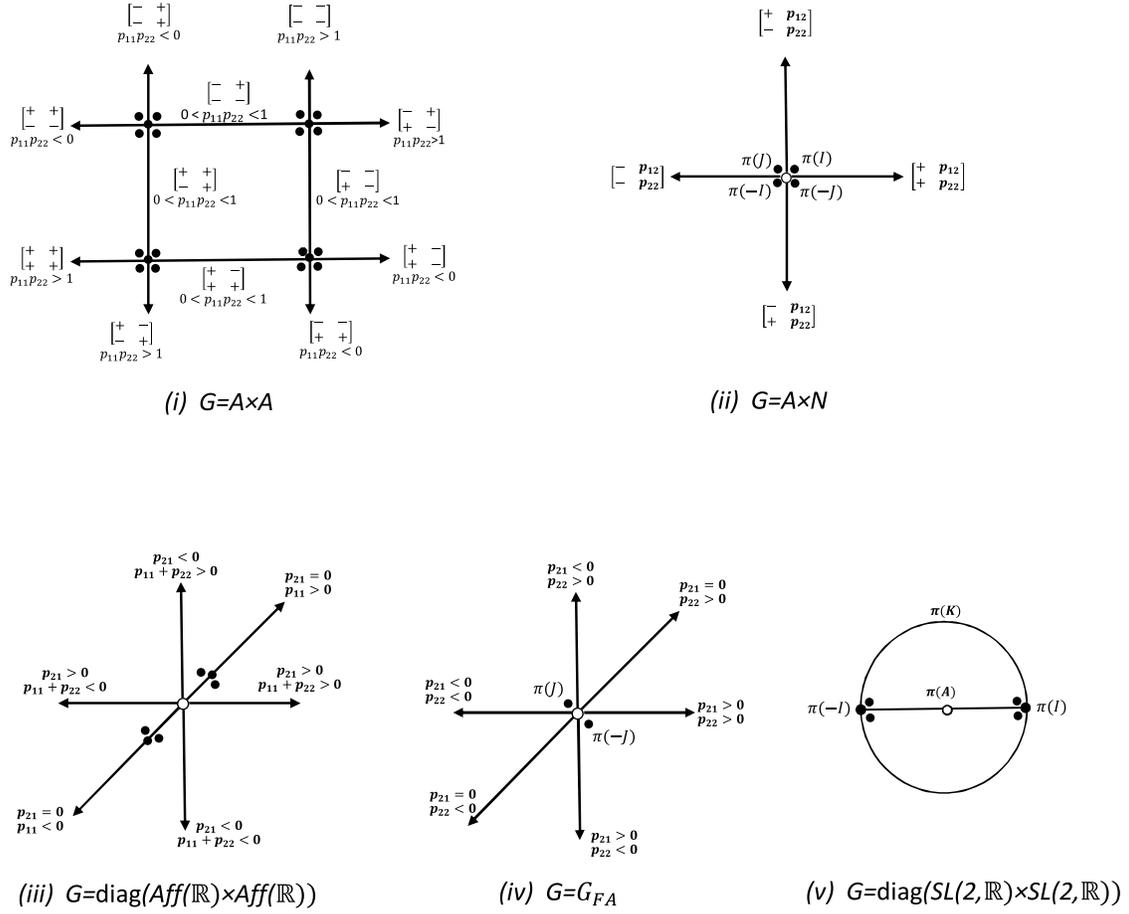}‎
‎\caption{The orbit space for the action of the five groups. In none of the above spaces, the topology is induced from $\R^2$. The signs
+ or - instead of $p_{ij}$ are used to denote the open subset of all
$p\in M$ where their $ij$-entry are positive or negative,
respectively.}\label{figAA}‎ ‎\end{figure}

 Let a connected and closed Lie subgroup $G$ of $Iso({\bf AdS_3})$ act isometrically and by \co on ${\bf AdS_3}$. Let $M={\bf AdS_3}$, and 
 $\pi : M\rightarrow M/G$ be the projection map on to the orbit space. Let $P$ and $S$ denote the set of principal and singular points in $M$, respectively.
 When the action is proper, the orbit space is Hausdorff and it is homeomorphic
 to one of the spaces $\R,\ S^1$ or $[0,+\infty)$, by Theorem \ref{proper}.
 In the case that the action is not proper, the orbit space may not be Hausdorff,
 and so the study becomes much more difficult. By Theorem \ref{1}-(iii) we know all the connected and closed Lie subgroups of the isometry group  acting by \co on $M$, up to conjugacy. So we consider the actions of these groups on $M$ to determine the orbit space up to homeomorphism as follows. To make it more clear we give presentations of five of these orbit spaces in Figure \ref{figAA}, up to homeomorphism, where their topologies are described in the context. Presenting the orbit spaces for other cases is more simple.

${\bf (i)}$ {\it The orbit space for the action of the Lie group} ${\bf G=A\times A .}$  Let $p$ and $p'$ be any two points of $M$. Then
\begin{equation}\label{orbit AA}G(p)=\{p_{11}e^{t-s}E_{11}+p_{12}e^{t+s}E_{12}+p_{21}e^{-(t+s)}E_{21}+p_{22}e^{-(t-s)}E_{22}|t,s\in\R\}.\end{equation}
By using (\ref{orbit AA}) and the proof of Theorem \ref{3.4.1}-(i) we have the following facts.

{\it Fact 1.} Let all of the entries of $p$ and $p'$ are nonzero.
Then $\pi(p)=\pi(p')$ if and only if $p_{11}p_{22}=p'_{11}p'_{22}$
and $p_{ij}p'_{ij}>0$, for each $i,j\in \{ 1,2\}$.

{\it Fact 2.} If $\pi(p)=\pi(p')$ then $p_{ij}p'_{ij}\geqslant0$. Furthermore, if $\pi(p)=\pi(p')$ and $p_{ij}=0$, for some $i,j\in \{ 1,2\}$, then $p'_{ij}=0$.

{\it Fact 3.} The four singular orbits are the orbits of $\pm I$ and
$\pm J$, where $J=E_{12}-E_{21}$. (The set of images of these orbits
in $M/G$, which we denote it by $\pi(S)$, is shown as the
intersection points in Figure \ref{figAA}).

{\it Fact 4.} A point $p$, where just one of its entries is zero, belongs to a space-like principal orbit.

{\it Fact 5.} If $p_{11}p_{22}>0$ or $p_{11}p_{22}<0$ then $G(p)$ is a Lorentzian surface. If $0<p_{11}p_{22}<1$ then $G(p)$ is a space-like surface.

By using these facts we determine the orbit space up to
homeomorphism. Facts 1 and 2 give a classification of points of the
orbit space. Let $P^{lo}$ and $P^{sp}$ denote the set of points (with no zero
entry) where their orbits are Lorentzian and space-like, respectively. Then Fact 5 implies that the sets
$\pi(P^{lo})$ (the eight unbounded open intervals emanated from the
intersection points in Figure \ref{figAA}-(i)) and $\pi(P^{sp})$ (the
four open intervals bounded by the intersection points in Figure
\ref{figAA}-(i)) are open in $M/G$ and the induced topology on them is
Hausdorff. Facts 2 and 4 shows that there are exactly sixteen
disjoint space-like orbits corresponding to the points $\pm I\pm
E_{12}$, $\pm I\pm E_{21}$, $\pm J\pm E_{11}$ and $\pm J\pm E_{22}$.
Denote the image of these orbits in the orbit space by
$\pi(P^{{\overline{sp}}})$, the sixteen points which are around the intersection points
in Figure \ref{figAA}-(i).

Consider the sequence $(g_n)=((A_n,A_n))$ in $G$. Then the four
sequences $g_n.(I\pm E_{21})$ and $g_{-n}.(I\pm E_{12})$, from
different space-like orbits, all converge to $I$, when $n\rightarrow
+\infty$. Hence $$\pi(I)\in \overline{\{\pi(I+
E_{12})\}}\cap\overline{\{\pi(I-E_{12})\}}\cap\overline{\{\pi(I+E_{21})\}}\cap\overline{\{\pi(I-
E_{21})\}}.$$ A similar discussion shows that each other point of $\pi(S)$
belongs to the intersection of the closures
of four one point sets of $\pi(P^{{\overline{sp}}})$. Hence any
neighborhood of each point of $\pi(S)$ contains four points of
$\pi(P^{{\overline{sp}}})$, two open intervals of $\pi(P^{lo})$ and
two open intervals of $\pi(P^{sp})$. An immediate consequence is
that the orbit space is neither Hausdorff nor locally Euclidean. We should point that the
sixteen points in Figure \ref{figAA}-(i) are not isolated points, since
any neighborhood of each of them contains an interval of
$\pi(P^{lo})$ and an interval of $\pi(P^{sp})$ (it may contain no
point of $\pi(S)$).

‎

${\bf (ii)}$ {\it The orbit space for the action of the Lie group} ${\bf G=A\times N.}$ By Theorem \ref{3.4.1}-(ii) and its proof
we know that there are just four distinct degenerate principal
orbits corresponding to the points $\pm I$ and $\pm J$, and each
other orbit $G(p)$, for which $p_{11}p_{21}\neq0$, is a Lorentzian
surface. In the case that $G(p)$ is a Lorentzian orbit, $p'\in G(p)$
if and only if $p_{11}p'_{11}>0$, $p_{21}p'_{21}>0$ and
$p_{11}p_{21}=p'_{11}p'_{21}$. This classifies the points in the
orbit space, and so we may consider the orbit space as the space in
Figure \ref{figAA}-(ii) with the following topology. The topology on the
four intervals emanated from the origin (which correspond to the set
of the images of Lorentzian orbits in $M/G$) is the subspace
topology induced from $\R^2$. Any neighborhood of each of the
remaining four points contains an interval on the horizontal axis
and an interval on the vertical axis, both emanated from the origin.
Hence the orbit space is not Hausdorff, but it is locally Euclidean.\\

${\bf (iii)}$ {\it The orbit space for the action of the Lie group} ${\bf G=N\times N.}$ By the proof of Theorem \ref{3.4.1}-(iii), the orbit
$G(p)$ is principal if and only if $p_{21}\neq0$. Hence $\pi(P)$ is
an open subset of $M/G$, and any neighborhood of each point of
$\pi(S)$ intersects $\pi(P)$. Let $G(q)$ be a singular orbit. Then
$\pi(q)=\pi(q')$ if and only if $q_{11}=q'_{11}$. If $G(p)$ is a
principal orbit, then $\pi(p)=\pi(p')$ if and only if
$p_{21}=p'_{21}$. Thus we may consider the orbit space $M/G$ as the
union of $x$ and $y$ axes in the plane from which the origin is removed, the $x$-axis for $\pi(P)$ and the $y$-axis for $\pi(S)$, with the following
topology. The topology on the $x$-axis
(without the origin) is the subspace topology induced from $\R^2$.
Any neighborhood of each point $\pi(q)$ of the $y$-axis contains
three disjoint intervals, one on the $y$-axis containing $\pi(q)$
and two on the $x$-axis emanating from the origin in opposite
directions. An immediate consequence is that the orbit space is neither
Hausdorff nor locally Euclidean.\\

${\bf (iv)}$ {\it The orbit space for the action of the Lie group} ${\bf G={\rm diag}(\Af\times \Af).}$ By the proof of Theorem
\ref{3.4.2}-(i) the orbit $G(p)$ is principal if and only if
$p\notin Aff(\R)$. If $G(p)$ is a
principal orbit, then $p'\in G(p)$ if and only if $p_{21}p'_{21}>0$
and $p_{11}+p_{22}=p'_{11}+p'_{22}$. This classifies the points of the orbit space in such a way that we may consider a homeomorphism from $\pi(P)$ to the union of $x$ and $y$ axes in the plane, from which the origin has been removed (see Figure \ref{figAA}-(iii). On the other hand, we may consider $\pi(S)$  as the union of the bisector of the first and third quarter of the plane (corresponding to the orbits which have been stated in Case (a) and (c) of the proof of
Theorem \ref{3.4.2}-(i)), and the four points (corresponding to the
four orbits stated in Case (b) of the proof of Theorem
\ref{3.4.2}-(i)) in Figure \ref{figAA}-(iii), with the topology which
is described as follows. Let $\pi(q_{11}E_{11}+q_{22}E_{22})$ be any
point in the bisector. If $q_{11}\neq q_{22}$, then any neighborhood
of $\pi(q_{11}E_{11}+q_{22}E_{22})$ contains three open intervals,
one on the bisector containing the point, and two on the $x$ and $y$
axes containing the points  $\pi(q_{11}E_{11}\pm
E_{21}+q_{22}E_{22})$. (note that $g_n.(q_{11}E_{11}\pm
E_{21}+q_{22}E_{22})\rightarrow q_{11}E_{11}+p_{22}E_{22}$, where
$(g_n)=((A_n,A_n))$). If $q_{11}=q_{22}$, i.e. $q\in \{\pm I\}$,
then any neighborhood of $\pi(q)$ contains three intervals as above,
and two near points of the four specified points (note that by Case (b) of the
proof of Theorem \ref{3.4.2}-(i), the four points correspond to
$\pi(\pm I\pm E_{12})$, and $g_{-n}.(I\pm E_{12})\rightarrow I$ and
$g_{-n}.(-I\pm E_{12})\rightarrow -I$). Any neighborhood of each of the four specified points contains two intervals, one on the $x$-axis and the other on $y$-axis, emanating from the origin. As a consequence, the
topology on the orbit space is neither Hausdorff nor locally
Euclidean.\\

${\bf (v)}$ {\it The orbit space for the action of the Lie group} ${\bf G=G_{FN}.}$ By Theorem \ref{3.4.2}-(ii) and its proof,
the action is free and there is no singular orbit. An orbit $G(p)$
is Lorentzian if $p_{21}\neq 0$ and it is degenerate if $p_{21}=0$.
There are two degenerate orbits corresponding to the points $\pm I$.
Hence the image of Lorentzian orbits is open and the induced
topology from $M/G$ is Hausdorff. For a Lorentzian orbit $G(p)$, by
a simple computation one gets that $p'\in G(p)$ if and only if
$p_{21}e^{-p_{22}/p_{21}}=p'_{21}e^{-p'_{22}/p'_{21}}$. Thus the
image of the Lorentzian orbits in $M/G$ is homeomorphic to $x$-axis in the plane from which the origin is removed. Imagine two points in above and below the origin corresponding to $\pi(I)$ and $\pi(-I)$. Then any neighborhood of $\pi(I)$ and $\pi(-I)$
intersects the $x$-axis in two open intervals, emanated from the origin
(the orbit space is known as an axis with two origins). Thus the
orbit space is not Hausdorff but it is locally Euclidean.\\

${\bf (vi)}$ {\it The orbit space for the action of the Lie group} ${\bf G=G_{FA}.}$ The following statements come from the proof
of Theorem \ref{3.4.2}-(iii). An orbit $G(p)$ is principal if and
only if $p\notin Aff(\R)$. Let $G(p)$ is a principal orbit. Then
$p'\in G(p)$ if and only if $p_{21}p'_{21}>0$ and $p_{22}=p'_{22}$.
Furthermore, $G(p)$ is a degenerate surface if and only if
$p_{21}\neq 0$ and $p_{22}=0$. This implies that there are only two
distinct degenerate principal orbits corresponding to the points
$\pm J$, and the image of other principal orbits in $M/G$ is
homeomorphic to the union of $x$ and $y$ axes where the origin has
been removed, with the induced topology from $\R^2$ (see Figure
\ref{figAA}-(iv)). Any $G$-invariant neighborhood of $J$ (resp. $-J$)
contains points where their $22$-entries can be negative, zero or
positive and their $21$-entry is negative (resp. positive). So any
neighborhood of $\pi(J)$ (rep. of $\pi(-J)$) contains two open
intervals emanating from the origin and containing $\pi(J\pm
\varepsilon E_{22})$ (resp. $\pi(-J\pm \varepsilon E_{22}$) for small
$\varepsilon>0$, one on the $x$-axis and the other on the $y$-axis in Figure
\ref{figAA}-(iv). Let $G(q)$ be a singular orbit. Then $q'\in G(q)$
if and only if $q_{21}=0$ and $q_{22}=q'_{22}$, which implies that there is a bijection from $\pi(S)$ to a line from which one of its points has been removed. Hence we may
consider $\pi(S)$, the image of singular points, as the bisector of
the first and third quarter of the plane in Figure \ref{figAA}-(iv), with the following topology. Since for any singular orbit $G(q)$ there
exists an element $g$ in $G$ such that $(g.q)_{12}=0$, so we may assume that
$q=q_{11}E_{11}+q_{22}E_{22}$. Hence any $G$-invariant neighborhood
of $q$ contains points with negative, zero and positive
$21$-entries. Thus any neighborhood of $\pi(q)$ contains three
disjoint open intervals, one on the bisector containing the point,
and two on the $x$ and $y$ axes containing the points $\pi(q\pm
E_{21})$. Note that $q\pm E_{21}$ are principal points and
$g_n.(q\pm E_{21})\rightarrow q$, where $g_n=(A_n,A_n)$. Hence
$\pi(q)\in \overline{(\pi(q+E_{21})}\cap \overline{(\pi(q-E_{21})}$. Thus the orbit space is neither Hausdorff nor locally Euclidean.\\

${\bf (vii)}$ {\it The orbit space for the action of the Lie group} ${\bf G=\Af\times A.}$ By Theorem \ref{3.4.3}-(i) the orbit
space consists of eight points, and by its proof the orbits of $\pm
I$ and $\pm J$ are the exceptional orbits, and the orbits of $\pm
I\pm E_{21}$ are the open orbits. Obviously, the subset $\{\pi(\pm
I\pm E_{21})\}\subset M/G$ has discrete topology.  Consider the
sequence $(g_n)=((A_n,A_n))$ in $G$. Then $g_n.(I\pm
E_{21})\rightarrow I$, and so $\pi(I)\in
\overline{\pi(I+E_{21})}\cap \overline{\pi(I-E_{21})}$. By a similar
discussion one sees that $\pi(-I)\in \overline{\pi(-I+E_{21})}\cap
\overline{\pi(-I-E_{21})}$, $\pi(J)\in \overline{\pi(I-E_{21})}\cap
\overline{\pi(-I-E_{21})}$ and $\pi(-J)\in
\overline{\pi(I+E_{21})}\cap \overline{\pi(-I+E_{21})}$. Since
$I\notin \overline{G(-I\pm E_{21})}$, $-I\notin \overline{G(I\pm
E_{21})}$, $J\notin \overline{G(\pm I+ E_{21})}$ and $-J\notin
\overline{G(I\pm E_{21})}$, so the following set is a basis for the
topology on $M/G$.
\begin{eqnarray*}\mathcal{B}&=&  \{  \{\pi(I+E_{21})\},\{\pi(I-E_{21})\},\{\pi(-I+E_{21})\},\{\pi(-I-E_{21})\},\\ & &\ \ \{\pi(I),\pi(I+E_{21}),\pi(I-E_{21})\},
\{\pi(-I),\pi(-I+E_{21}),\pi(-I-E_{21})\}, \\ & & \ \
\{\pi(J),\pi(I-E_{21}),\pi(-I-E_{21})\},\{\pi(-J),\pi(I+E_{21}),\pi(-I+E_{21})\}
\}.\end{eqnarray*} This shows that the orbit space is not either
Hausdorff or locally Euclidean.\\

${\bf (viii)}$ {\it The orbit space for the action of the Lie group} ${\bf G=\Af \times N.}$ By Theorem \ref{3.4.3}-(ii) the
orbit space consists of four points, and by its proof the orbits of
$\pm I$ are the exceptional orbits, and the orbits of $I\pm E_{21}$
are the open orbits. Clearly, the subset $\{\pi(I\pm
E_{21})\}\subset M/G$ has discrete topology. If
$(g_n)=((F_{n,1-e^n},N_{1-e^n}))$ and
$(g'_n)=((F_{n,-1-e^n},N_{-1-e^n}))$, then
$g_n.(I+E_{21})\rightarrow I$ and $g'_n.(I+E_{21})\rightarrow -I$.
If $(h_n)=((F_{n,e^n-1},N_{e^n-1}))$ and
$(h'_n)=((F_{n,e^n+1},N_{e^n+1}))$ then $h_n.(I-E_{21})\rightarrow
I$ and $h'_n.(I-E_{21})\rightarrow -I$. Hence
$\{\pi(I),\pi(-I)\}\subset
\overline{\pi(I+E_{21})}\cap\overline{\pi(I-E_{21})}$. Thus the
following set is a basis for the topology on $M/G$.
\begin{eqnarray*}\mathcal{B}&=&  \{  \{\pi(I+E_{21})\},\{\pi(I-E_{21})\},\\ & &\ \ \{\pi(I),\pi(I+E_{21}),\pi(I-E_{21})\},
\{\pi(-I),\pi(I+E_{21}),\pi(I-E_{21})\} \}.\end{eqnarray*} As a
consequence, the orbit space is not either Hausdorff of locally
Euclidean.\\

${\bf (ix)}$ {\it The orbit space for the action of the Lie group} ${\bf G=\Af \times \Af.}$ It is easily seen that the orbit
space is homeomorphic to that of the case $G=\Af \times N$.\\

${\bf (x)}$ {\it The orbit space for the action of the Lie group} ${\bf G=diag(SL(2,\R)\times SL(2,\R)).}$ In this case we know
from the proof of Theorem \ref{3.4.2}-(iii) that $p'\in G(p)$ if and
only if $p$ and $p'$ are conjugate within $SL(2,\R)$. Let $p\notin
\{\pm I\}$. Then $p$ is conjugate to some element of $N\cup(-N)$,
$A\cup (-A)$ or $K$.

Let $p\in A\cup (-A)-\{\pm I\}$. Then $p'$ is conjugate to $p$ if
and only if $p'=p^{-1}$ or $p'=p$. The set of
hyperbolic elements is an open subset of $M$. Hence each of the $\pi(A)-\{I\}$
and $\pi(-A)-\{-I\}$ is homeomorphic to an open interval in $\R^2$. So we may consider
them as the open diameter of the circle from which the center is removed (see Figure \ref{figAA}-(v)).

Let $p\in K-\{\pm I\}$. Then $p$ is conjugate to itself only. Hence
$\pi(K)=S^1$.

Any two elements $p$ and $p'$ of $N-\{I\}$ (and of $-N-\{-I\}$) are
conjugate if and only if $p_{12}p'_{12}>0$. This implies that
$\pi(N)=\{\pi(I),\pi(I+E_{12}),\pi(I-E_{12})\}$ and
$\pi(-N)=\{\pi(-I),\pi(-I+E_{12}),\pi(-I-E_{12})\}$. Since any
$G$-invariant neighborhood of each parabolic element intersects the
sets of hyperbolic and elliptic elements in open sets, so any
neighborhood of each of the specified four points in Figure
\ref{figAA}-(v) contains its near intersection point and two
intervals of $\pi(A-\{I\})$ (or of $\pi(-A-\{-I\})$) and of
$\pi(K)-\{\pm I\}$, both emanating from the intersection point.

It is easily seen that $\pi(I)\in
\overline{\pi(I+E_{12})}\cap\overline{\pi(I-E_{12})}$ and
$\pi(-I)\in \overline{\pi(-I+E_{12})}\cap\overline{\pi(-I-E_{12})}$.
Hence each neighborhood of $\pi(I)$ (resp. of $\pi(-I)$) contains
the union of one open interval, two open arcs, all emanating from
the point as in Figure \ref{figAA}-(v), and two points $\pi(I\pm
E_{12})$ (resp. $\pi(-I\pm E_{21})$). As a consequence, the orbit
space is neither Hausdorff nor locally Euclidean. But the induced
topology on $M/G-\{\pi(\pm I),\pi(\pm I\pm E_{12})\}$ is Hausdorff
and locally Euclidean.

${\bf (xi)}$ {\it The orbit space for the action of the Lie group} ${\bf G=G_{FF}.}$ By Theorem \ref{3.4.2}-(iii) and its proof
the following statements hold. There are two distinct open orbits
and uncountably many singular orbits, where each of them is
diffeomorphic to $\R$. An orbit $G(q)$ is a singular orbit if and
only if $q_{21}=0$. Let $q$ be a typical singular point. We may
assume that $q_{12}=0$, since there is a point $q'$ in $G(q)$ where
$q'_{12}=0$. Then the orbits of the points $q\pm E_{21}$ are
the open orbits and $g_n.(q\pm E_{21})\rightarrow q$, where
$(g_n)=((A_n,A_n))$. On the other hand, we have $q'\in G(q)$ if and
only if $q'_{21}=0$ and $q_{11}=q'_{11}$.  So we may consider the
image of singular orbits, say $\pi(S)$, as a line where one of its points has been removed. Call this point origin, and consider the image of two open orbits in $M/G$, say $\pi(I\pm E_{21})$, as two points in above and below the origin, with the following topology. Two one
point sets $\{\pi(I+E_{21})\}$ and $\{\pi(I-E_{21})\}$ are open in
$M/G$. Any neighborhood of each point $\pi(q)$ of the line contains
an open interval around the point and the two points $\{\pi(I\pm
E_{21})\}$. Clearly, the orbit space is neither Hausdorff nor locally Euclidean. (note that this
space is not an axis with two origins, since its topology is
different).

Thus we have the following proposition as a consequence of this section. 

\begin{pro}
Let ${\bf adS_3}$ be of cohomogeneity one under the action of a connected
and closed Lie subgroup $G\subset Iso({\bf adS_3})$. If the action is not
proper, then the following statements hold.

(i) The orbit space is not Hausdorff.

(ii) The orbit space is locally Euclidean if and only if $G$ is conjugate to either $A\times N$ or $G_{FN}$.

(iii) The orbit space is a finite space if and only if $G$ is conjugate to one of the groups $\Af \times A$, $\Af \times N$ or $\Af \times \Af$.
\end{pro}

{\bf Acknowledgement:} The author should gratefully thank the referee for invaluable suggestions
  leading to the improvement of the paper.

\bibliographystyle{amsplain}

\begin{thebibliography}{10}
 \bibitem{Adams}    S. Adams, Dynamics on Lorentz manifolds, {\it world Scientific}, 2001.
\bibitem{Alee} D.V. Alekseevsky, On a proper action of a Lie group,
{\it Uspekhi Mat. Nauk}  {\bf 34} (1979), 219-220.
\bibitem{AA} A.V. Alekseevsky and D.V. Alekseevsky, $G$-manifolds
 with one dimensional orbit space,{\it Adv. Sov. Math.} {\bf{8}} (1992)
 1-31.
  \bibitem{AK1} P. Ahmadi and S.M.B. Kashani, Cohomogeneity one de
 Sitter space $S^n_1$, {\it Acta Math. Sin.} Vol.{\bf 26} No. {\bf
 10 } (2010) 1915-1926.
 \bibitem{H^3_1} P. Ahmadi and S.M.B. Kashani, Cohomogeneity one
 anti de-Sitter space $H^3_1$, {\it Bull. Iran. Math. Soc.}, Vol. {\bf 35},  No.
 {\bf 1} (2009) 221-233.
\bibitem{Ber}  L. Berard-Bergery, Sur de nouvells vari$\acute{e}$t$\acute{e}$
       riemanniennes d'Einstein, {\it Inst. $\acute{E}$lie~Cartan} No.{\bf 6} (1982), 1-60.
\bibitem{Br}  G.E. Bredon, Introduction to compact
 transformation groups. {\it Academic Press, New York,} 1972.
  \bibitem{DK}   J.J. Duistermaat, J.A.C. Kolk, Lie Groups,
 {\it Springer}, 2000.
\bibitem{He}  S. Helgason, Differential geometry, Lie groups and
 symmetric spaces, {\it Academic Press, Inc.} 1978.
 \bibitem{K}  R.S. Kulkarni, Proper actions and
 pseudo-Riemannian space forms, {\it Adv. in Math.} {\bf 40}
 (1981), 10-51.
\bibitem{Kn}  A.W. Knapp, Lie groups beyond an introduction,
second edition, {\it Progress in Math.} Volume 140, 2002.
\bibitem{Mos}  P.S. Mostert, On a compact Lie group acting on a manifold,
         {\it Ann. Math.} Vol.{\bf{65}}, No.{\bf{3}} (1957), 447-455.
\bibitem{PT}   R.S. Palais and CH.L. Terng, A general theory of
canonical forms, {\it Trans. Am. Math. Soc}. {\bf{300}} (1987),
771-789.
\bibitem{PT2}   R.S. Palais and CH.L. Terng, Critical Point Theory and Submanifold Geometry, Lectture Notes in Mathematics, {\it Springer-Verlag}, 1988.
\bibitem{PS}  F. Podesta and A. Spiro: Some topological properties of chomogeneity
one manifolds with negative curvature,{\it Ann. Global Anal. Geom.}
{\bf 14} (1996), 69-79.
\bibitem{S}  C. Searle, Cohomogeneity and positive curvature in
low dimension, {\it Math. Z.} {\bf 214} (1993), 491-498.
\bibitem{V}  L. Verdiani, Invariant
       metrics on \co manifolds, {\it Geom. Ded.} {\bf 77} (1999) 77-111.



 \end{thebibliography}

\end{document}